\documentclass{amsart} 
\title{Counting colored planar maps free-probabilistically}
\author{Abdelmalek Abdesselam}
\email{malek@virginia.edu}
\address{University of Virginia Department of Mathematics,
P. O. Box 400137, Charlottesville, VA 22904-4137, USA}
\author{Greg W.\ Anderson}
\email{gwanders@umn.edu}
\address{University of Minnesota, Minneapolis, MN 55455, USA}
\thanks{A.A. was supported in part by the National Science Foundation under grant DMS 0907198}
\date{August 10, 2012}
\newcommand{\coU}{{\mathrm{coU}}}
\newcommand{\Pos}{{\mathrm{Pos}}}
\newcommand{\HHH}{{\mathcal{H}}}

\newcommand{\dd}{{\mathrm{d}}}

\newcommand{\Sym}{{\mathrm{Sym}}}

\newcommand{\Match}{{\mathrm{Match}}}

\newcommand{\Color}{{\mathrm{Color}}}
\newcommand{\Arb}{{\mathrm{Arb}}}
\newcommand{\zbold}{{\mathbf{z}}}
\newcommand{\kbold}{{\mathrm{k}}}
\newcommand{\jbold}{{\mathbf{j}}}
\newcommand{\Poly}{{\mathrm{Poly}}}

\newcommand{\Splice}{{\mathrm{Splice}}}

\newcommand{\Ibold}{{\mathbf{I}}}

\newcommand{\Cov}{{\mathrm{Cov}}}

\newcommand{\Zbold}{{\mathbf{Z}}}
\newcommand{\weight}{{\mathrm{wt}}}

\newcommand{\T}{{\mathrm{T}}}       
 
\newcommand{\Xbold}{{\mathbf{X}}}

\newcommand{\norm}[1]{{\left\Vert#1\right\Vert}}
\newcommand{\Perm}{{\mathrm{Perm}}}
\newcommand{\Mat}{{\mathrm{Mat}}}
\newcommand{\AAA}{{\mathcal{A}}}

\newcommand{\trace}{{\mathrm{tr}}}
\newcommand{\Map}{{\mathrm{Map}}}

\newcommand{\RR}{{\mathbb{R}}}
\newcommand{\st}{{\textup{\mbox{\scriptsize s.t.}}}}
\newcommand{\ii}{{\mathrm{i}}}
\newcommand{\Ebold}{{\mathbf{E}}}
\newcommand{\Part}{{\mathrm{Part}}}
\newcommand{\CC}{{\mathbb{C}}}
\def\bbone{{\mathchoice {\rm 1\mskip-4mu l} {\rm 1\mskip-4mu l}
{\rm 1\mskip-4.5mu l} {\rm 1\mskip-5mu l}}}
\newcommand{\zero}{{\mathbf{0}}}
\newcommand{\one}{{\mathbf{1}}}
\newcommand{\ZZ}{{\mathbb{Z}}}
\newcommand{\ibold}{{\mathbf{i}}}
\newcommand{\Tfrak}{{\mathfrak{T}}}

\newcommand{\Ffrak}{{\mathfrak{F}}}
\newcommand{\Gfrak}{{\mathfrak{G}}}
\newtheorem{Lemma}[subsubsection]{Lemma}
\newtheorem{Theorem}[subsubsection]{Theorem}
\newtheorem{TheoremBis}[subsection]{Theorem}
\newtheorem{Proposition}[subsubsection]{Proposition}
\newtheorem{Corollary}[subsubsection]{Corollary}
\newtheorem{CorollaryBis}[subsection]{Corollary}
\theoremstyle{remark}
\newtheorem{Remark}[subsubsection]{Remark}
\newtheorem{Example}[subsubsection]{Example}

\begin{document}

\begin{abstract}
Our main result is an explicit operator-theoretic formula for the number of
colored planar maps with a fixed set of stars each of which has a fixed set of half-edges
with fixed coloration.
The formula bounds the number of such colored planar maps well enough to prove convergence near the origin of generating functions arising  naturally
in the matrix model context. Such convergence is known but the proof of convergence
proceeding by way of our main result is relatively simple. Besides Voiculescu's  generalization of Wigner's semicircle law,
our main technical tool is an integration identity
representing the joint cumulant of several functions of a Gaussian random vector.
The latter identity in the case of cumulants of order $2$
reduces to one well-known as a means to prove the Poincar\'{e} inequality.
We derive the identity by combining the heat equation with the so-called BKAR formula
from constructive quantum field theory and rigorous statistical mechanics.
\end{abstract}

\maketitle

\section{Introduction}
Let $\Perm(n)$ denote the group of bijective maps of $\langle n\rangle=\{1,\dots,n\}$ to itself. Let $c(\theta)$ denote the number of cycles into which $\theta\in \Perm(n)$ 
decomposes.
Let  $\Match(n)\subset\Perm(n)$ denote the subset consisting of fixed-point-free elements
of order $2$.
Let $\Color(n)$ denote the set consisting of maps of $\langle n\rangle$ to itself. 
Given any $\theta\in \Perm(n)$ and $\gamma\in \Color(n)$, we define finite sets
\begin{eqnarray}\label{equation:MapDef}
\Map(\theta,\gamma)&=&
\left\{\iota\in \Match(n)\left\vert
\begin{array}{l}
\mbox{$\theta$ and $\iota$ generate a group}\\
\mbox{acting transitively on $\langle n\rangle$}\\
\mbox{and furthermore $\gamma\circ \iota=\gamma$}
\end{array}\right.\right\}\;\;\mbox{and}\\
\label{equation:Map0Def}
\Map_0(\theta,\gamma)&=&\{\iota\in \Map(\theta,\gamma)\mid c(\theta)+c(\theta\iota)=2+n/2\}.
\end{eqnarray}
Our main result, namely Theorem \ref{Theorem:BigKahuna} below,
gives a novel free-probabilistic representation of the cardinality $|\Map_0(\theta,\gamma)|$.
We will locate our main result more precisely with respect to the literature 
after formulating it in the next section.

The significance of the quantity $|\Map_0(\theta,\gamma)|$ in random
matrix theory is explained by the following well-known formula which we owe to the physicists.
We refer to   \cite{Zvonkin} for background on matrix integrals and maps.
See also \cite{GuionnetMaurelSegala} for a good explanation to mathematicians of the physicists' viewpoint on matrix models. 
Recall that an $N$-by-$N$ {\em GUE matrix} $\Xi$
is a random hermitian matrix whose entries have a centered Gaussian joint distribution satisfying $\Ebold\Xi(i,j)\Xi(i',j')=\delta_{ij'}\delta_{i'j}$.

\begin{TheoremBis}['t Hooft \cite{tHooft}]\label{Theorem:tHooft}
Fix $\theta\in \Perm(n)$ and $\gamma\in \Color(n)$. 
Let
$$\theta=(i_{1,1}\cdots i_{1,n_1})\cdots (i_{k,1}\cdots i_{k,n_k})\;\;\;(n_1+\cdots+n_k=n)$$
be the decomposition of $\theta$ into cycles.
Let $\Xi^{(N)}_1,\dots,\Xi^{(N)}_n$ be independent $N$-by-$N$ GUE matrices.
Then we have
\begin{equation}\label{equation:MatricialFormOftHooft}
|\Map_0(\theta,\gamma)|
=\lim_{N\rightarrow\infty}
\frac{\kappa\left(\trace (\Xi^{(N)}_{\gamma(i_{1,1})}\cdots \Xi^{(N)}_{\gamma(i_{1,n_1})}),
\dots, \trace (\Xi^{(N)}_{\gamma(i_{k,1})}\cdots \Xi^{(N)}_{\gamma(i_{k,n_k})})\right)}
{\displaystyle N^{2+n/2-k}},
\end{equation}
where 
$\kappa(\cdot)$ is the joint cumulant functional.\end{TheoremBis}
\noindent Recall that for $\CC$-valued random variables $X_1,\dots,X_k$ with absolute moments of all
orders the {\em joint cumulant} $\kappa(X_1,\dots,X_k)\in \CC$ is defined by the relation
\begin{eqnarray}\label{equation:JointCumulantDefPhys}
&&\log \sum_{\nu_1,\dots,\nu_k=0}^\infty \Ebold(X_1^{\nu_1}\cdots X_k^{\nu_k})
\frac{t_1^{\nu_1}\cdots t_k^{\nu_k}}{\nu_1!\cdots \nu_k!}\\
\nonumber&=&\sum_{\nu_1,\dots,\nu_k=1}^\infty
\kappa\left(\underbrace{X_1,\dots,X_1}_{\nu_1},\dots,
\underbrace{X_k,\dots,X_k}_{\nu_k}\right)\frac{t_1^{\nu_1}\cdots t_k^{\nu_k}}{\nu_1!\cdots \nu_k!}
\end{eqnarray}
standing between formal power series. For $\kappa(\cdot)$ one also has an expression
(see equation \eqref{equation:JointCumulantDef} below)
involving the M\"{o}bius function of the lattice of set partitions. 

Of course Theorem \ref{Theorem:tHooft}
is only a faint reflection of the full matrix model picture.
Formula \eqref{equation:MatricialFormOftHooft} merely peels off the leading term of the $1/N$-expansion. We will not be considering the higher order terms in this paper---but see Remark \ref{Remark:Susceptible} below.

Physicists have a visually appealing interpretation of $\iota\in \Map_0(\theta,\gamma)$
which we recall to justify the title of the paper.
One draws a diagram as follows.
First one marks down some vertices in the plane indexed by the $\theta$-cycles.
Then out of each vertex one draws half-edges indexed by the elements of the
corresponding 
$\theta$-cycle,  arranging them in the circular order dictated by $\theta$
and ``coloring'' them by $\gamma$, thus forming a ``star.''
Finally, one joins each half-edge to its same-colored mate via the perfect matching $\iota$
to form a whole edge, taking care that connecting paths do not cross---the numerical condition $c(\theta)+c(\theta\iota)=2+n/2$ 
guarantees the possibility of drawing the picture in the plane without crossings. 
The resulting graph embedded in the plane (up to some abuse of language)  is a {\em colored planar map}.
See  \cite{LandoZvonkin} for  background on graphs embedded in surfaces.

Actually our work in this paper consists entirely of an analysis of the right side of equation \eqref{equation:MatricialFormOftHooft}. Thus the reader could 
without loss of comprehension of our arguments
 take \eqref{equation:MatricialFormOftHooft} as the definition of the quantity $|\Map_0(\theta,\gamma)|$.

Although we must defer the statement of Theorem \ref{Theorem:BigKahuna},
we can immediately state a couple of corollaries to it in order to convey its flavor.
We use the notation $\bbone\{\cdots\}$ for the characteristic function of the condition
between braces. For constant $\gamma\in \Color(n)$ (the ``monochrome'' case)
we  just write $\Map_0(\theta)$ instead of $\Map_0(\theta,\gamma)$.

\begin{CorollaryBis}\label{Corollary:BigKahuna}
Let $n_1,\dots,n_k>0$ be integers. Put $n=\sum n_i$ and $p=\prod n_i$.
Let $\theta\in \Perm(n)$ have cycles of length $n_1,\dots,n_k$.
We have\begin{equation}\label{equation:BigKahunaBis}
|\Map_0(\theta)|
\leq pn^{k-2} 2^{n-2k+2}\bbone\{n\geq 2k-2\}.
\end{equation}
\end{CorollaryBis}
\noindent This estimate may be new. In any case the 
method of proof is surely new.
See \S\ref{subsubsection:KahunaBisProof} below for proof of the corollary,
and see Proposition \ref{Proposition:ArbBound} immediately following
for a more precise version of this bound.

\begin{CorollaryBis}\label{Corollary:PositiveRadius}
Let $n_1,\dots,n_k>0$ be integers. For  integers $\nu_1,\dots,\nu_k>0$
let $\theta_{\nu_1,\dots,\nu_k}\in \Perm(\sum \nu_in_i)$ have cycles of length
$\underbrace{n_1,\dots, n_1}_{\nu_1}, \dots,
\underbrace{n_r,\dots, n_k}_{\nu_k}$.
Then the generating function
\begin{equation}\label{equation:TheGeneratingFunction}
\sum_{\nu_1=1}^\infty\cdots
\sum_{\nu_k=1}^\infty\left|\Map_0\left(\theta_{\nu_1,\dots,\nu_k}\right)\right|
\frac{ z_1^{\nu_1}\cdots z_k^{\nu_k}}{\nu_1!\cdots \nu_k!}
\end{equation}
converges for $(z_1,\dots,z_k)$ in a neighborhood of the origin in $\CC^k$.
\end{CorollaryBis}
\noindent This result is certainly not new.
It follows e.g. from results of \cite{ErcolaniMcLaughlin}
or \cite{GuionnetMaurelSegala}.
But the relatively elementary character of our proof is novel. \proof 
By Corollary \ref{Corollary:BigKahuna} the series in question is majorized by the series
$$\sum_{\nu_1=1}^\infty\cdots\sum_{\nu_k=1}^\infty \;
\frac{(\sum_{i=1}^k \nu_i)!}{\nu_1!\cdots \nu_k!}\;\;
\frac{\left(\sum_{i=1}^k \nu_i n_i\right)^{\sum_{i=1}^k \nu_i} }{(\sum_{i=1}^k\nu_i)!}\;
\prod_{i=1}^k(n_i 2^{n_i}z_i)^{\nu_i}
$$
and the latter clearly has a positive radius of convergence. 
\qed

Now we turn to a topic which seems at first glance only mildly relevant.
Let $\zeta,\zeta^{(1)},\zeta^{(2)}\in \RR^n$ be independent random vectors, each with i.i.d.\ standard Gaussian entries. 
Let $f,g:\RR^n\rightarrow\CC$ be adequately nice
functions. The identity
\begin{equation}\label{equation:UrArboreal}
\Cov(f(\zeta),g(\zeta))=\int_0^1\Ebold\left(\nabla f(\zeta^{(1)})\cdot \nabla g\left(t\zeta^{(1)}+\sqrt{1-t^2}\,\zeta^{(2)}\right)\right)\,\dd t
\end{equation}
is well-known. Notably, \eqref{equation:UrArboreal} implies the Poincar\'{e} inequality with the best constant.
For  background, further applications of \eqref{equation:UrArboreal} and more references,
see e.g. \cite{BobkovGotzeHoudre}.
The main technical result of the paper, namely Theorem \ref{Theorem:MalliavinApp} below,
is a generalization of \eqref{equation:UrArboreal} holding for joint cumulants
of arbitrary order.
We will locate our main technical result more precisely with respect to the literature after formulating it in the next section. We will derive Theorem \ref{Theorem:MalliavinApp}  from the Brydges-Kennedy-Abdesselam-Rivasseau (BKAR) formula \cite{BrydgesK,AbdesselamR1}
which has appeared  in the context of constructive quantum field theory and rigorous statistical mechanics.

Here is the plan of the paper.
In \S\ref{section:JointCumulantFunctionsGauss} we 
formulate  our main result, Theorem~\ref{Theorem:BigKahuna}.
We also state our main technical result, Theorem \ref{Theorem:MalliavinApp}.
In \S\ref{section:BKAR},
we review the BKAR formula (see Theorem \ref{Theorem:BKAR} below), give pointers to the related mathematical physics literature, and finally we use BKAR and the heat equation to prove Theorem \ref{Theorem:MalliavinApp}.
In \S\ref{section:OperatorTheoreticCharacterization}
we prove Theorem \ref{Theorem:BigKahuna} by combining
Theorem \ref{Theorem:MalliavinApp}, Voiculescu's multi-matrix generalization of Wigner's semicircle law \cite[Thm. 2.2]{Voiculescu} 
and Theorem \ref{Theorem:tHooft}. 

\section{Formulation and discussion of the  main results}
\label{section:JointCumulantFunctionsGauss}

\subsection{Co-ultrametrics, forests and random matrices}
We introduce notions needed to state our main result, Theorem \ref{Theorem:BigKahuna} below,
and our main technical result, Theorem \ref{Theorem:MalliavinApp} below.
As in the introduction, we write $\langle k\rangle=\{1,\dots,k\}$.

\subsubsection{Co-ultrametrics}
Let $A$ be a $k$-by-$k$ real symmetric matrix whose entries $A(i,j)$ satisfy the following
conditions for all $i_1,i_2,i_3\in \langle k\rangle$:
\begin{eqnarray}
\label{equation:Co-ultrametric2}
A(i_1,i_1)&=&1,\\\label{equation:Co-ultrametric4}
A(i_1,i_2)&\in&[0,1],\\\label{equation:Co-ultrametric3}
A(i_1,i_3)&\geq &\min(A(i_1,i_2),A(i_2,i_3)).
\end{eqnarray}
We call $A$ a $k$-by-$k$ {\em co-ultrametric}. We choose this terminology because \eqref{equation:Co-ultrametric3} is the reversal of the ultrametric inequality familiar (say) to number theorists. Let $\coU(k)$ denote the set of $k$-by-$k$ co-ultrametrics. 
\begin{Lemma}\label{Lemma:NoLoops}
Fix $A\in \coU(k)$.  Let $j_1,\dots,j_r\in \langle k\rangle$ be a sequence with $r\geq 2$.
Then the two smallest numbers on the list $A(j_1,j_2),\dots,A(j_r,j_1)$ are equal.
\end{Lemma}
\noindent This is the ``co'' version of a standard fact about ultrametrics. (``Every triangle is isosceles.'')
\proof On the one hand, after cyclically permuting the indices $j_1,\dots,j_r$, we may assume that
$A(j_1,j_r)=A(j_r,j_1)\leq \min(A(j_1,j_2),\dots,A(j_{r-1},j_r))$.
On the other hand, the reverse inequality holds by repeated application of \eqref{equation:Co-ultrametric3}.
\qed

\begin{Lemma}\label{Lemma:PosDef}
Every $A\in \coU(k)$ is positive semidefinite,
and moreover positive definite
if all off-diagonal entries of $A$ are less than $1$.
\end{Lemma}
\proof For every $A\in \coU(k)$ and $t\in [0,1)$ there exists unique $A^{(t)}\in \coU(k)$
with entries $A^{(t)}(i,j)=\bbone\{A(i,j)>t\}$. One checks immediately that $A^{(t)}$ is
the indicator of the graph of an equivalence relation in $\langle k\rangle$ and in particular
is positive semidefinite. Furthermore, if the off-diagonal entries of $A$ are all less than one,
then $A^{(t)}$ is the identity matrix for $t$ near $1$.
The trivial formula $A=\int_0^1 A^{(t)}\,dt$ completes the proof. \qed

\subsubsection{Gapless co-ultrametrics}
We say that $A\in \coU(k)$ has {\em articulation}
equal to the cardinality of the set $\{A(i,j)\mid i,j=1,\dots,k\}\setminus \{0,1\}$.
We say that $A\in \coU(k)$ has {\em co-articulation}
equal to the number of equivalence classes in the set $\langle k\rangle$ for the equivalence
relation $i\sim j\Leftrightarrow A(i,j)>0$. We say that $A\in \coU(k)$ is {\em gapless}
if articulation and co-articulation sum to $k$. 

\subsubsection{Graphs (especially trees)}
In this paper a {\em graph} $\Gfrak=(V,E)$ is a pair consisting of (i) a finite  set $V$ of {\em vertices}
and (ii) a set $E$ of (unoriented) {\em edges} each of which is a two-element subset of $V$. 
We say that $\Gfrak=(V,E)$ is a {\em tree} if $\Gfrak$ is connected and $|E|=|V|-1$. 
We say that a tree $\Tfrak=(V,E)$ {\em spans} a set $S$ if $S=V$.
Recall the generating function identity
\begin{equation}\label{equation:Kirchhoff}
\sum_{\begin{subarray}{c}\Tfrak:\,\mbox{\scriptsize tree}\\
\mbox{\scriptsize spanning $\langle k\rangle$}
\end{subarray}}\;\;\;
\prod_{\mbox{\scriptsize $\{i,j\}$: edge  of $\Tfrak$}}
x_ix_j\;\;\;\;=\;\;\;\;(x_1+\cdots+x_k)^{k-2}x_1\cdots x_k
\end{equation}
which one obtains by specializing Kirchhoff's  matrix-tree theorem. 
In particular, there are exactly $k^{k-2}$ trees spanning the set $\langle k\rangle$.
The graphs of interest in this paper are mostly trees and if not trees then {\em forests}, i.e.,
disjoint unions of trees, or equivalently, graphs with no circuits.
When we say that $\Ffrak$ is a forest in $\langle k\rangle$ we mean that $\Ffrak=(V,E)$
is a forest with vertex set $V=\langle k\rangle$.

In fact the concept of gaplessness is closely related to forests,
as follows.

\begin{Proposition}\label{Proposition:Obvious}
Let $\Ffrak=(V,E)$ be forest in $\langle k\rangle$.
Let $\{x_e\}_{e\in E}$ be a family of distinct numbers selected from the open unit interval $(0,1)$.
Then there exists unique $A\in \coU(k)$ such that \begin{eqnarray}\label{equation:Co-ultrametric5}
&&\mbox{$A(i,j)=x_{\{i,j\}}$ for every edge $\{i,j\}\in E$ and}\\\label{equation:Co-ultrametric6}
&&\mbox{$A$ is gapless of articulation $|E|$.}
 \end{eqnarray}
\end{Proposition}
\proof  By Lemma \ref{Lemma:NoLoops},
for distinct $i,j\in \langle k\rangle$ connected by some walk in $\Ffrak$,
any co-ultrametric $A$ satisfying \eqref{equation:Co-ultrametric5} and \eqref{equation:Co-ultrametric6} also satisfies
\begin{equation}\label{equation:Co-ultrametric7}
A(i,j)=\min\left(x_{\{j_1,j_2\}},\dots,x_{\{j_{r-1},j_r\}}\right)
\end{equation}
where $j_1\cdots j_r$ is the unique geodesic walk in $\Ffrak$ from $i$ to $j$.
(It is at this point that the assumption of distinctness of the numbers $x_e$ enters crucially.)
Now consider the following two equivalence relations $\sim_1$ and $\sim_2$ in $\langle k\rangle$.
Let $i\sim_1 j$ if and only if $i$ and $j$ are joined by a walk in $\Ffrak$.
Let $i\sim_2 j$ if and only if $A(i,j)>0$. We have $i\sim_1 j\Rightarrow i\sim_2 j$.
But furthermore, $\sim_1$ has $k-|E|$ equivalence classes because $\Ffrak$ is a forest
and $\sim_2$ has $k-|E|$ equivalence classes by \eqref{equation:Co-ultrametric6}.
Thus the equivalence relations $\sim_1$ and $\sim_2$ coincide. It follows that
\begin{equation}\label{equation:Co-ultrametric8}
\mbox{$A(i,j)=0$ if $i,j\in \langle k\rangle$ are not joined by a walk in $\Ffrak$.}
\end{equation}
Thus uniqueness is settled. Conditions \eqref{equation:Co-ultrametric2}, 
 \eqref{equation:Co-ultrametric7} 
and \eqref{equation:Co-ultrametric8} define a $k$-by-$k$ real symmetric matrix which one can easily verify does indeed belong to $\coU(k)$ and satisfy \eqref{equation:Co-ultrametric5} and \eqref{equation:Co-ultrametric6}.
Thus existence is settled. \qed

The preceding result has an easily derived converse.
\begin{Proposition}\label{Proposition:ObviousConverse}
Let $A\in \coU(k)$ be gapless of articulation $\ell$.
Let 
$$\{i_1,j_1\},\dots,\{i_\ell,j_\ell\}$$ be two-element subsets of $\langle k\rangle$
such that 
$$0<A(i_1,j_1)<\dots<A(i_\ell,j_\ell)<1.$$
Then the graph $(\langle k\rangle,\{\{i_1,j_1\},\dots,\{i_\ell,j_\ell\}\})$ is a forest in $\langle k\rangle$.
\end{Proposition}
\proof Lemma \ref{Lemma:NoLoops} rules out the possibility of circuits.
\qed

\subsubsection{The weight matrix attached to a forest}\label{section:weights}
Let $\Ffrak=(V,E)$ be a forest in $\langle k\rangle$.
Let $\{x_e\}_{e\in E}$ be i.i.d. random variables uniformly distributed in $(0,1)$.
We define \linebreak $\weight_\Ffrak\in \coU(k)$ to be the random matrix which with probability $1$ is uniquely determined  by
Proposition \ref{Proposition:Obvious}. By Lemma \ref{Lemma:PosDef},
with probability $1$, the matrix
$\weight_\Ffrak$ is positive definite and thus has a unique positive definite square root $\sqrt{\weight_\Ffrak}$.
We introduce the random matrices $\weight_\Ffrak$ just to have a convenient compact notation with which to handle integrals on various pieces of the space $\coU(k)$.

\subsection{Spanning trees and splicing involutions}
We introduce some partly group-theoretical
and partly arboreal notions needed to formulate Theorem~\ref{Theorem:BigKahuna} below.
\subsubsection{Permutations}
Recall from the introduction that $\Perm(n)$ denotes the group of bijective maps of $\langle n\rangle=\{1,\dots,n\}$ to itself
and $c(\theta)$ denotes the number of cycles into which $\theta\in \Perm(n)$ decomposes.
We call $\tau\in \Perm(n)$ an {\em involution} if $\tau^2=1$.

\begin{Proposition}[Riemann-Hurwitz bound]
\label{Proposition:RiemannHurwitzBound}
Let $\sigma_0,\sigma_1,\sigma_\infty\in \Perm(n)$  generate a group of permutations
acting transitively on $\langle n\rangle$
and satisfy  $\sigma_0\sigma_1\sigma_\infty=1$.
Then we have $c(\sigma_0)+c(\sigma_1)+c(\sigma_\infty)\leq n+2$.
\end{Proposition}
\proof The triple $(\sigma_0,\sigma_1,\sigma_\infty)$ is an example of a {\em constellation}.
See \cite{LandoZvonkin} for background on constellations. 
To the triple $(\sigma_0,\sigma_1,\sigma_\infty)$ one naturally attaches
an \linebreak $n$-sheeted covering of the Riemann sphere branched at $0$, $1$ and $\infty$.
The genus $g$ of that covering, i.e., the number of its ``handles,'' satisfies the {\em Riemann-Hurwitz formula}
$2g-2=-2n+\sum_{i=0,1,\infty} (n-c(\sigma_i))$.
The result follows because $g\geq 0$. \qed

\begin{Remark} Note that the numerical condition figuring
in the definition \eqref{equation:Map0Def} of $\Map_0(\theta,\gamma)$
touches the bound enunciated in Proposition \ref{Proposition:RiemannHurwitzBound}
and thus implies ``planarity'' of elements of $\Map_0(\theta,\gamma)$.
\end{Remark}

\subsubsection{Splicing involutions}
Let $\nu:\langle n\rangle\rightarrow\langle k\rangle$ be an onto function.
Let $\Tfrak$ be a tree spanning $\langle k\rangle$.
We define the set of {\em splicing involutions}
$$\Splice_\Tfrak(\nu)\subset \Perm(n)$$ 
indexed by $\Tfrak$ and $\nu$
to be the subset consisting of
involutions $\tau$  admitting a factorization 
$$\tau=(i_1,j_1)\cdots(i_{k-1}j_{k-1})$$
into disjoint transpositions such that
$$\Tfrak=(\langle k\rangle,\{\{\nu(i_1),\nu(j_1)\},\dots,\{\nu(i_{k-1}),\nu(j_{k-1})\}\}).$$

\begin{Proposition}\label{Proposition:SplicingCyclicity}
Let $\theta\in \Perm(n)$ be any permutation and put $k=c(\theta)$.
Let \linebreak $\nu:\langle n\rangle\rightarrow \langle k\rangle$ be any $\theta$-invariant onto map.
Let $\Tfrak$ be any tree spanning $\langle k\rangle$.
Then for all $\tau\in \Splice_\Tfrak(\nu)$ the composite permutation $\theta\tau$ is cyclic.
\end{Proposition}
\proof Since onto, the map $\nu$ in effect enumerates the $\theta$-cycles. It follows in turn
that the disjoint transpositions into which $\tau$ factors provide enough linkages
between the $\theta$-cycles to make the Cayley graph of the pair $(\theta,\tau)$ connected.
Thus $\theta$ and $\tau$ generate a group of permutations
acting transitively on $\langle n\rangle$. 
By Proposition \ref{Proposition:RiemannHurwitzBound}
it follows that $c(\theta\tau)=1$.
\qed

\begin{Proposition}\label{Proposition:SpliceCount}
Let $\nu:\langle n\rangle\rightarrow\langle k\rangle$ be onto. We have
\begin{equation}\label{equation:LocalSpliceCountBis}
\sum_{\begin{subarray}{c}
\Tfrak:\;\textup{\mbox{\scriptsize tree}}\\
\textup{\mbox{\scriptsize spanning $\langle k\rangle$}}
\end{subarray}}|\Splice_\Tfrak(\nu)|=\left\{\begin{array}{rl}
\displaystyle \frac{ (n-k)! \prod_{i=1}^k n_i} { (n-2k+2)!}&\mbox{if $n\geq 2k-2$,}\\
0&\mbox{if $n<2k-2$,}
\end{array}\right.
\end{equation}
where $n_i=|\nu^{-1}(i)|$ for $i=1,\dots,k$. 
\end{Proposition}
\proof For each tree $\Tfrak$ spanning $\langle k\rangle$,
let $\vec{\Tfrak}$ be the set of ordered pairs \linebreak $(i,i')\in \langle k\rangle^2$
such that $\{i,i'\}$ is an edge of $\Tfrak$, let $\pi_\Tfrak=((i,i')\mapsto i):\vec{\Tfrak}\rightarrow\langle k\rangle$,
and put $\mu_{\Tfrak,i}=|\pi_\Tfrak^{-1}(i)|$ for $i=1,\dots,k$.
Now the set $\Splice_\Tfrak(\nu)$ is in evident bijective
correspondence with the set of one-to-one maps $\psi:\vec{\Tfrak}\rightarrow\langle n\rangle$ such that $\nu\circ \psi=\pi_\Tfrak$. Thus we have
$$
(\mbox{LHS of \eqref{equation:LocalSpliceCountBis}})=\sum_{\begin{subarray}{c}
\Tfrak:\;\mbox{\scriptsize tree}\\
\mbox{\scriptsize spanning $\langle k\rangle$}
\end{subarray}}\prod_{i=1}^k \frac{n_i!}{(n_i-\mu_{\Tfrak,i})!},
$$
where it is understood that we set $\frac{n_i!}{(n_i-\mu_{\Tfrak,i})!}$ equal to zero
if $\mu_{\Tfrak,i}>n_i$.
Note furthermore that we can rewrite identity \eqref{equation:Kirchhoff}
in the form
$$\sum_{\begin{subarray}{c}\Tfrak:\,\textup{\mbox{\scriptsize tree}}\\
\textup{\mbox{\scriptsize spanning}}\;
\langle k\rangle
\end{subarray}} \;\;\prod_{\vec{e}\in \vec{\Tfrak}}x_{\pi_\Tfrak(\vec{e})}
=\sum_{\begin{subarray}{c}\Tfrak:\,\textup{\mbox{\scriptsize tree}}\\
\textup{\mbox{\scriptsize spanning}}\;
\langle k\rangle
\end{subarray}}\;\;\prod_{i=1}^k x_i^{\mu_{\Tfrak,i}}
=(x_1+\cdots +x_k)^{k-2}x_1\cdots x_k.
$$
Finally, we have
$$\left(\frac{\partial}{\partial x_1}+\cdots +\frac{\partial}{\partial x_k}\right)^{k-2}\frac{\partial^k}
{\partial x_1\cdots \partial x_k}(x_1^{n_1}\cdots x_k^{n_k})\bigg\vert_{x_1=\cdots=x_k=1}=
(\mbox{RHS of \eqref{equation:LocalSpliceCountBis}}),$$
which finishes the proof.
\qed

\subsubsection{Colored splicing involutions}
Let $\nu:\langle n\rangle\rightarrow\langle k\rangle$ be an onto function.
Let $\Tfrak$ be a tree spanning $\langle k\rangle$.
Given also $\gamma\in \Color(n)$ we put
$$\Splice_\Tfrak(\nu,\gamma)=\{\tau\in \Splice_\Tfrak(\nu)\mid \gamma\circ \tau=\gamma\},$$
thus defining the colored version of $\Splice_\Tfrak(\nu)$.

\subsection{Statement of the main result}

\subsubsection{Data} 
Fix $\theta\in \Perm(n)$ and $\gamma\in \Color(n)$.
Our aim is to formulate a result representing the quantity $|\Map_0(\theta,\gamma)|$ free-probabilistically. For this purpose we put $k=c(\theta)$ to abbreviate
and we fix a $\theta$-invariant onto map $\nu:\langle n\rangle\rightarrow\langle k\rangle$.
\subsubsection{The canonical splicing polynomial}
Let $\CC\langle\{\Xbold_i\}_{i=1}^n\rangle$ be the noncommutative polynomial algebra
generated by a family $\{\Xbold_i\}_{i=1}^n$ of $n$ independent noncommuting variables.
Let $\Tfrak$ be any tree spanning $\langle k\rangle$.
We then define
\begin{eqnarray*}
&&\Poly_{\theta,\gamma,\nu,\Tfrak}\\
&=&\sum_{\tau\in \Splice_\Tfrak(\nu,\gamma)}
\sum_{\begin{subarray}{c}
i_1,\dots,i_n\in \langle n\rangle\\
\st \;\theta\tau=(i_1\cdots i_n)\\
\mbox{\scriptsize and}\;i_1=1
\end{subarray}}
\Xbold_{i_1}^{\bbone\{i_1=\tau(i_1)\}}
\cdots \Xbold_{i_n}^{\bbone\{i_n=\tau(i_n)\}}\in\CC\langle\{\Xbold_i\}_{i=1}^n\rangle.
\end{eqnarray*}
The inner sum has exactly one term 
by Proposition \ref{Proposition:SplicingCyclicity}, so that the number of monomials
in $\Poly_{\theta,\gamma,\nu,\Tfrak}$ equals $|\Splice_\Tfrak(\nu,\gamma)|$.
Clearly, $\Poly_{\theta,\gamma,\nu,\Tfrak}$  is homogeneous of degree $n-2k+2$ in $\Xbold_1,\dots,\Xbold_n$
and in particular vanishes identically unless \linebreak $n\geq 2k-2$.
One should think of $\Poly_{\theta,\gamma,\nu,\Tfrak}$ as a sort of noncommutative generating function
representing the set $\Splice_\Tfrak(\nu,\gamma)$.

\subsubsection{A transformation of the canonical splicing polynomial}
Let 
$\{\{\Zbold(i,j)\}_{i=1}^k\}_{j=1}^n$ be another family of independent noncommutative algebraic variables.
Given a tree $\Tfrak$ spanning $\langle k\rangle$ put
$$\Zbold_\Tfrak(i,j)=\sum_{i'=1}^k \sqrt{\weight_\Tfrak}(i,i')\Zbold(i',j)\;\;\mbox{for $(i,j)\in \langle k\rangle\!\times\!\langle n\rangle$,}
$$
(which is random) and in turn put
$$
\Arb_{\theta,\gamma,\nu}=\sum_{
\begin{subarray}{c}
\mbox{\scriptsize trees $\Tfrak$}\\
\mbox{\scriptsize spanning $\langle k\rangle$}
\end{subarray}}
\Ebold \left(\Poly_{\theta,\gamma,\nu,\Tfrak}\left|_{\begin{subarray}{l}
\Xbold_i=\Zbold_\Tfrak(\nu(i),\gamma(i))\\
\mbox{\scriptsize for $i\in \langle n\rangle$}
\end{subarray}}\right)\right.\in \CC\langle\{\{\Zbold(i,j)\}_{i=1}^k\}_{j=1}^n\rangle,
$$
where the expectation is computed term-by-term, i.e., one first expands the integrand as a
sum of finitely many monomials in the variables $\Zbold(i,j)$ with bounded random coefficients,
and then takes the expectation of each of the finitely many nonzero coefficients.
We remark that $\Arb_{\theta,\gamma,\nu}$ actually is independent of $\nu$.
This fact will not be needed in the sequel.

\subsubsection{Free standard semicircular variables}
Let $(\AAA,\varphi)$ be a faithful tracial \linebreak $C^*$-probability space.
Let 
$\{\{\zbold(i,j)\}_{i=1}^k\}_{j=1}^n$ be a family of free 
standard semicircular variables in $(\AAA,\varphi)$. 
A quick introduction to free probability theory more than adequate for our purposes here
can be found in \cite[Chap.\ 5]{AGZ}. 
For discussion at  length and in depth see \cite{NicaSpeicher} or \cite{VDN}.

We are ready to state the main result of the paper.
\begin{Theorem}\label{Theorem:BigKahuna}
Notation and assumptions are as above.  We have
\begin{equation}\label{equation:BigKahuna}
|\Map_0(\theta,\gamma)|
=\varphi \left(\Arb_{\theta,\gamma,\nu}\left|_{\begin{subarray}{l}\Zbold(i,j)=\zbold(i,j)\\
\mbox{\scriptsize for $(i,j)\in \langle k\rangle\!\times\!\langle n\rangle$}
\end{subarray}}\right).\right.
\end{equation}
\end{Theorem}
\noindent 
The proof of Theorem \ref{Theorem:BigKahuna} 
will be completed in \S\ref{section:OperatorTheoreticCharacterization} below.
\subsubsection{Proof of Corollary \ref{Corollary:BigKahuna}}\label{subsubsection:KahunaBisProof}
The corollary is entirely subsumed by the following more precise technical result.

\begin{Proposition}\label{Proposition:ArbBound} 
In the setup for Theorem \ref{Theorem:BigKahuna} and with
 $\norm{\cdot}$ denoting the norm on the $C^*$-algebra $\AAA$, we have
\begin{eqnarray}\label{equation:ArbBound}
&&\norm{\Arb_{\theta,\gamma,\nu}\left|_{\begin{subarray}{l}\Zbold(i,j)=\zbold(i,j)\\
\mbox{\scriptsize for $(i,j)\in \langle k\rangle\!\times\!\langle n\rangle$}
\end{subarray}}\right.}\\
\nonumber&\leq &\left\{\begin{array}{rl}
\displaystyle 2^{n-2k+2}\frac{ (n-k)! \prod_{i=1}^k n_i} { (n-2k+2)!}&\mbox{if $n\geq 2k-2$,}\\
0&\mbox{if $n<2k-2$.}
\end{array}\right.
\end{eqnarray}
\end{Proposition}
\proof Given any 
tree $\Tfrak$ spanning $\langle k\rangle$ we define
$$\zbold_\Tfrak(i,j)=\sum_{i'=1}^k \sqrt{\weight_\Tfrak}(i,i')\,\zbold(i',j)\in \AAA\;\;
\mbox{for $i\in \langle k\rangle$ and $j\in \langle n\rangle$}
$$
(which is random).
With probability one, because the diagonal entries of the matrix $\weight_\Tfrak$ are identically equal to $1$,
for each fixed $i\in \langle k\rangle$ 
the family $\{\zbold_\Tfrak(i,j)\}_{j=1}^n$ 
is free standard semicircular and in particular $\norm{\zbold_\Tfrak(i,j)}=2$ for all $i$ and $j$.
Consequently we have
$$\norm{\Poly_{\theta,\gamma,\nu,\Tfrak}\bigg\vert_{\begin{subarray}{l}
\Xbold_i=\zbold_\Tfrak(\nu(i),\gamma(i))\\
\mbox{\scriptsize for}\;i\in \langle n\rangle
\end{subarray}}}
\leq 2^{n-2k+2}|\Splice_\Tfrak(\nu,\gamma)|,
$$
almost surely.
Trivially, for a tree $\Tfrak$ spanning $\langle k\rangle$ 
and any bounded linear functional $\psi\in \AAA^*$ we have the formula
$$k^{k-2}\Ebold\psi\left(\Poly_{\theta,\gamma,\nu,\Tfrak}\left|_{\begin{subarray}{l}
\Xbold_i=\zbold_\Tfrak(\nu(i),\gamma(i))\\
\mbox{\scriptsize for}\;i\in \langle n\rangle
\end{subarray}}\right)\right.=\psi \left(\Arb_{\theta,\gamma,\nu}\left|_{\begin{subarray}{l}\Zbold(i,j)=\zbold(i,j)\\
\mbox{\scriptsize for $(i,j)\in \langle k\rangle\!\times\!\langle n\rangle$}
\end{subarray}}\right).\right.
$$
By Hahn-Banach we can choose $\psi$ of norm $1$ to make the right side above equal to the left side of \eqref{equation:ArbBound}.
We then get the result by Proposition \ref{Proposition:SpliceCount}.
\qed

\begin{Remark} If $k=1$ and $\theta=(i_1\cdots i_n)$,
formula \eqref{equation:BigKahuna} reduces to the statement
\begin{equation}\label{equation:MatricialFormOftHooftBis}
|\Map_0(\theta,\gamma)|=\varphi(\zbold(1,\gamma(i_1))\cdots\zbold(1,\gamma(i_n))),
\end{equation}
which nowadays, free probability theory taken for granted, one interprets as an instance of the semicircular analogue of  the Wick formula.
\end{Remark}
\begin{Remark} In a continuation of the previous remark, we note
that when specialized to the case $k=1$,
Theorems \ref{Theorem:tHooft} and \ref{Theorem:BigKahuna} in tandem are equivalent to  Voiculescu's result \cite[Thm. 2.2]{Voiculescu}.  \end{Remark}

\begin{Remark}\label{Remark:PreSpeicherBis} 
We turn now to the case $k=2$.
In this case Theorem \ref{Theorem:BigKahuna} provides a formula for the same quantity as does \cite[Thm.\ 5.3]{MingoSpeicher}.
The formulas look rather different.
On the one hand, \cite[Thm.\ 5.3]{MingoSpeicher} involves 
an extension of free probability theory to the second order (cyclic Fock spaces and annular pair partitions intervene).
On the other hand, Theorem \ref{Theorem:BigKahuna} is phrased  in terms of the usual first order theory,
albeit with certain operations of differentiation and integration added on.
 Detailed relations between the two types of formulas remain to be worked out.
\end{Remark}

\begin{Remark}\label{Remark:SpeicherBis}
In the case $k\geq 2$
there is some overlap, at least thematically, between our work here
and that in \cite{MingoSpeicherBis}. We note, for example, that both theories have in common a key role for permutations subject to what we call here
Riemann-Hurwitz bounds. But the program of \cite{MingoSpeicherBis} is much more ambitious, seeking to understand limiting behavior of
joint cumulants of traces of monomials in independent unitary invariant random matrices (not just GUE matrices), and it does so
by enlarging the foundations of free probability to all orders.
More precise connections to the less sophisticated theory worked out here remain to be elucidated.
\end{Remark}

\begin{Remark}
One can construct noncommutative
({\em a priori} possibly signed) measures
through consideration
of generating functions of the form \eqref{equation:TheGeneratingFunction}
for small values of the parameters.
The formidable analyses of \cite{GuionnetMaurelSegala} 
and \cite{GuionnetShlyakhtenko} taken together show (among many other things) 
that under the evident necessary
conditions those measures are positive
and engender von Neumann algebras
with useful and interesting properties. An important motivation for our work was to find a relatively elementary approach to the positivity phenomenon independent of matrix models. We did not succeed in finding it but we nonetheless hope Theorem \ref{Theorem:BigKahuna}
could provide some clues in this direction.
\end{Remark}

\begin{Remark}
 Tutte developed over several decades
 and in many papers an approach to the enumeration of planar maps based on ``well-labeled trees''
and generating functions. For just one influential example, see \cite{Tutte}.
We also mention  \cite{BerBous}, \cite{BousJeh} and \cite{BousScha} as recent papers
in Tutte's line of descent. The relationship between the trees appearing here
and those of Tutte is apparently not very direct but deserves investigation.
\end{Remark}

\subsection{Formulation of the main technical result}
We state the main technical result of the paper, namely Theorem \ref{Theorem:MalliavinApp} below.

\subsubsection{The matrix space $\Mat_{k\times n}$}
Let $\Mat_{k\times n}$ denote the space of $k$-by-$n$ matrices with real entries.
We equip $\Mat_{k\times n}$ with the
inner product $(x,y)=x\cdot y=\trace\,( x^\T y)$ and we put $\norm{x}=(x,x)^{1/2}$.
We write $\Mat_k=\Mat_{k\times k}$.
Let $\Ibold_k\in \Mat_k$ denote the $k$-by-$k$ identity matrix.
Let $D_{ij}$ denote differentiation of smooth functions on $\Mat_{k\times n}$
with respect to the entry in position $(i,j)$.
Given a matrix $Z\in \Mat_{k\times n}$, we often
denote its entry in row $i$ and column $j$ by $Z(i,j)$.
We tend to choose notation of the latter type
when there are further indices to keep track of.

\subsubsection{Classes of functions defined on $\Mat_{k\times n}$}
We say that $f:\Mat_{k\times n}\rightarrow\CC$ is {\em smooth}
if infinitely differentiable; 
{\em polynomial} if expressible as a polynomial in the entries;
 {\em of exponential growth} if there exist constants $c_1$ and $c_2$
depending only on $f$
such that $|f(x)|\leq c_1e^{c_2\norm{x}}$ for all 
$x\in \Mat_{k\times n}$; and {\em tame} if smooth and furthermore
partial derivatives of all orders have exponential growth. We carry the preceding terminology
over to functions $f:\RR^n\rightarrow\CC$ via the
identification $\RR^n=\Mat_{1\times n}$.

\subsubsection{Tensor products}
Let $f_1,\dots,f_k:\RR^n\rightarrow\CC$ be tame functions.
The tame function $f:\Mat_{k\times n}\rightarrow\CC$ defined by the formula
$$
f\left(\left[\begin{array}{ccccc}
x_{11}&\cdots&x_{1n}\\
\vdots&&\vdots\\
x_{k1}&\cdots&x_{kn}
\end{array}\right]\right)=\prod_{i=1}^kf_i(x_{i1},\dots,x_{in})
$$
will be  denoted by $f_1\otimes \cdots \otimes f_k$.

\subsubsection{Differential operators indexed by trees}
For each tree $\Tfrak$ spanning the set $\langle k\rangle$
and smooth function $f:\Mat_{k\times n}\rightarrow\CC$
we define
$$L_\Tfrak f=\left(\prod_{
\{i,i'\}:\; \mbox{\scriptsize edge of $\Tfrak$}}\;\; \sum_{j=1}^n D_{ij}D_{i'j}\right)f.
$$

\begin{Theorem}\label{Theorem:MalliavinApp}
Let  $\zeta\in \RR^n$ and $Z\in \Mat_{k\times n}$  have i.i.d.\ standard Gaussian entries. 
Assume that for all trees $\Tfrak$ spanning $\langle k\rangle$ the random matrix $\weight_\Tfrak$ 
is independent of $Z$.
Then for any tame functions $f_1,\dots,f_k:\RR^n\rightarrow\CC$ 
we have
\begin{equation}\label{equation:MalliavinApp}
\kappa(f_1(\zeta),\dots,f_k(\zeta))=
\sum_{\begin{subarray}{c}
\textup{\mbox{\scriptsize trees $\Tfrak$}}\\
\textup{\mbox{\scriptsize spanning $\langle k\rangle$}}
\end{subarray}}\Ebold\,(L_\Tfrak(f_1\otimes \cdots \otimes f_k))(\sqrt{\weight_\Tfrak}\,Z),
\end{equation}
where $\kappa(\cdot)$ is the joint cumulant functional.
\end{Theorem}
\noindent 
The proof of Theorem \ref{Theorem:MalliavinApp} will be completed in \S\ref{subsection:ProofOfMalliavinApp} below after we have introduced the BKAR formula.
\begin{Example}
\label{Example:Case2} 
In the case $k=2$ of Theorem \ref{Theorem:MalliavinApp}, consider the unique
tree $\Tfrak$ spanning the set $\langle 2\rangle$.
Let $U$ be a random variable uniformly distributed in $(0,1)$ and independent of the matrix $Z$.
We then have
$$\weight_\Tfrak\stackrel{d}{=}\left[\begin{array}{cc}
 1&U\\
 U&1
 \end{array}\right]\;\;
 \mbox{and}\;\;
 \sqrt{\weight_\Tfrak}\,Z\stackrel{d}{=}\left[\begin{array}{cc}
 1&0\\
 U&\sqrt{1-U^2}
 \end{array}\right]Z.$$
Thus \eqref{equation:MalliavinApp} in the case $k=2$ reduces to formula \eqref{equation:UrArboreal}. 
\end{Example}

\begin{Remark}
We will prove Theorem \ref{Theorem:MalliavinApp} by exploiting
the so-called BKAR formula \cite{BrydgesK,AbdesselamR1}, thus placing it in the 
long tradition of ``arboreal'' representations of truncated correlation functions in statistical mechanics and quantum field theory.
But Theorem \ref{Theorem:MalliavinApp} is also in many respects similar, say, to \cite[Thm.\ 5.1]{NourdinPeccati}  and thus in principle also part of a long tradition of ``arboreal'' representations of cumulants of
Wiener processes.
To set up a ``dictionary'' linking these two distinct traditions using the specific example of Theorem \ref{Theorem:MalliavinApp} is an interesting problem, if largely one of exposition.
Perhaps the vantage point of the paper \cite{RotaStochastic} permits a unified view.
\end{Remark}

\begin{Remark}
Theorem \ref{Theorem:MalliavinApp} complements  the influential central limit theorem \linebreak \cite[Thm.\ 4.2]{Chatterjee} by supplying explicit means to compute limiting covariance  as well as higher corrections in typical random matrix applications. \end{Remark}

\begin{Remark}
On the one hand,
Theorem \ref{Theorem:MalliavinApp} uses more general test-functions than we need to make our application to counting colored planar maps.
For the latter purpose only polynomial test-functions are actually needed.
On the other hand, the regularity assumptions in Theorem \ref{Theorem:MalliavinApp}
are clearly far stronger than necessary. The minimal hypotheses concerning regularity
remain to be worked out. 
\end{Remark}

\begin{Remark}
It would be interesting to obtain a cumulant representation
for functions of a vector uniformly distributed on the discrete cube 
generalizing \linebreak \cite[Thm. 3.1]{BobkovGotzeHoudre} 
in the same sense that Theorem \ref{Theorem:MalliavinApp}
generalizes \eqref{equation:UrArboreal}.
\end{Remark}
\begin{Remark}\label{Remark:PrePreSpeicherBis}
It would be interesting to find an analogue of Theorem \ref{Theorem:MalliavinApp} with the Gaussian random vector replaced by a Haar-distributed element
of the $N$-by-$N$ unitary group.  This problem is suggested by Remark \ref{Remark:SpeicherBis} above.
\end{Remark}

\section{Application of the BKAR formula}\label{section:BKAR}
In this section we state the BKAR formula, derive  a ``connected'' version of the formula,
provide some background and references,
and finally prove Theorem \ref{Theorem:MalliavinApp}.

\subsection{Statement  of the BKAR formula}
The BKAR formula is a combinatorial amplification of the Fundamental Theorem of Calculus.
In its present form the identity first appeared in~\cite[Thm. III.1]{AbdesselamR1}, building on the previous work of
Brydges and Kennedy~\cite{BrydgesK}.   For the reader who is not familiar with the literature of constructive quantum field theory and rigorous statistical mechanics, the most accessible reference for BKAR  is the recent introduction with complete proofs given in~\cite{AbdesselamNote}. 
\subsubsection{Note on notation}
Notation used here is not exactly the same as in \cite{AbdesselamNote} but the differences are of a trivial nature. 
For example we employ $k$-by-$k$ symmetric matrices here
to play the role assigned in \cite{AbdesselamNote}  to arrays indexed by two-element subsets of $\langle k\rangle$. The reader should not have great difficulty translating.

\subsubsection{Notation related to symmetric matrices}
Let $\Sym_k$ denote the space of \linebreak $k$-by-$k$ matrices with real entries.
Let $\Sym_k^+\subset \Sym_k$ denote the closed set of positive semidefinite matrices
and let $\Pos_k\subset \Sym_k$ denote the open set of positive definite matrices.
For any smooth function $f:U\rightarrow\CC$ defined on an open subset $U\subset \Sym_k$,
and given matrices $X\in U$ and $Y\in \Sym_k$, we define the directional derivative by $\nabla_Yf(X)=\frac{d}{dt}f(X+tY)\big\vert_{t=0}$.

\subsubsection{Differential operators indexed by forests}
Let $\Ffrak=(V,E)$ be a forest in $\langle k\rangle$.
For any smooth function $f:U\rightarrow\CC$ defined on an open subset $U\subset \Sym_k$
we define
\begin{equation}
\dd_\Ffrak f=\left(\prod_{\{i,j\}\in E}\nabla_{e_{ij}+e_{ji}}\right)f
\end{equation}
where $e_{ij}\in \Mat_k$ denotes the $k$-by-$k$ elementary matrix with entry $1$ in position $(i,j)$
and $0$ elsewhere.
If $E=\emptyset$, then $\dd_\Ffrak$ is simply the identity operator.

\begin{Theorem}[BKAR formula]\label{Theorem:BKAR}
For any smooth function $f:\Sym_k\rightarrow\CC$,
\begin{equation}\label{equation:BKAR}
f\left(\left[\begin{array}{ccc}
1&\dots&1\\
\vdots&&\vdots\\
1&\dots&1\end{array}\right]\right)
=\sum_{\begin{subarray}{c}
\textup{\mbox{\scriptsize forests $\Ffrak$}}\\
\textup{\mbox{\scriptsize in $\langle k\rangle$}}
\end{subarray}}\Ebold(\dd_\Ffrak f)(\weight_\Ffrak).
\end{equation}
\end{Theorem}
\noindent 
Proofs of Theorem \ref{Theorem:BKAR} can be found in~\cite[Thm. III.1]{AbdesselamR1},
\cite[Thm. VIII.2]{BrydgesM} and \cite{AbdesselamNote}. 

\subsection{The connected BKAR formula}
To prove Theorem \ref{Theorem:BigKahuna} we need to rewrite equation \eqref{equation:BKAR} so as to make the right side into a sum over trees rather than forests. We undertake that exercise in this subsection
after recalling a few facts about set partitions and
introducing appropriate notation.  See for example \cite{Rota} or \cite{Shiryaev} for background on the topic of set partitions, M\"{o}bius functions and joint cumulants.

\subsubsection{The lattice $\Part(k)$}
Let $\Part(k)$ denote the lattice of set partitions of $\langle k\rangle$.
By definition the elements of $\Part(k)$ are disjoint families of nonempty subsets of $\langle k\rangle$
whose union is $\langle k\rangle$. Members of a set partition are called its {\em blocks}.
Given partitions $\Phi,\Psi\in \Part(k)$ we say that $\Phi$ is a {\em refinement} of $\Psi$
and write $\Phi\leq \Psi$ if and only if for all $A\in \Phi$ there exists $B\in \Psi$ such that $A\subset B$.
The set $\Part(k)$ is partially ordered by the relation of  refinement
and thus becomes a {\em lattice}. Let 
$$\one_k=\{\{1,\dots,k\}\}\;\;\mbox{and}\;\;\zero_k=\{\{1\},\dots,\{k\}\},$$
which are the maximal and minimal elements of $\Part(k)$, respectively, with respect to the refinement partial order.

\subsubsection{Matrix representation of partitions}
Given $\Phi\in \Part(k)$ we  define a matrix $[\Phi]\in \coU(k)$ by the rule
that $[\Phi](i,j)=\sum_{A\in \Phi}\bbone\{(i,j)\in A\times A\}$.
For example
$$[\{\{1\},\{2,3\}\}]=\left[\begin{array}{ccc}
1\\
&1&1\\
&1&1\end{array}\right].$$
Note that in general
$$[\one_k]=\left[\begin{array}{ccc}
1&\dots&1\\
\vdots&&\vdots\\
1&\dots&1\end{array}\right]
\;\;\mbox{and}\;\;
[\zero_k]=\left[\begin{array}{ccc}
1\\
&\ddots\\
&&1\end{array}\right]=\Ibold_k.$$
It is clear that  the map $(\Phi\mapsto [\Phi]):\Part(k)\rightarrow\coU(k)$ is one-to-one.
More precisely, for any partitions $\Phi,\Psi\in \Part(k)$ 
we have $\Phi\leq \Psi$ if and only if $[\Phi](i,j)\leq [\Psi](i,j)$ for $i,j=1,\dots,n$.

\begin{Remark}\label{Remark:CanBeShown} It can be shown that for every $A\in \coU(k)$ there exist 
\begin{itemize}
\item a unique
chain $\{\Phi_0>\cdots>\Phi_\ell\}$ in $\Part(k)$ and 
\item unique numbers $0<t_1<\cdots<t_\ell<1$
\end{itemize}
such that 
$$A=[\Phi_\ell]+\sum_{i=1}^\ell t_i([\Phi_{i-1}]-[\Phi_i]).$$
One proves this by elaborating the idea of the proof of Lemma \ref{Lemma:PosDef}. Thus
$\coU(k)$ has a canonical simplicial decomposition indexed by chains in $\Part(k)$
and in this qualified sense the points of $\coU(k)$ of the form $[\Phi]$ for $\Phi\in \Part(k)$
are its vertices. But note that $\coU(k)$ is not convex for $k\geq 3$ nor is it a topological manifold. It has in general a rather complicated and interesting shape. In the simplest nontrivial case $k=3$
it has the fish-like shape of  three triangular fins stuck onto a line segment backbone.
\end{Remark}

\subsubsection{The M\"{o}bius function of the lattice of set partitions}
\label{subsubsection:MoebiusSymbol}
The {\em M\"{o}bius function}  
$$\mu(\cdot:\cdot)=\mu_k(\cdot:\cdot):\Part(k)\times \Part(k)\rightarrow\ZZ$$
is by definition characterized by the {\em M\"{o}bius inversion formula}
\begin{equation}\label{equation:MoebiusInversion}
\sum_{\Theta\in \Part(S)}\bbone\{\Pi\leq \Theta\}\;\mu(\Theta:\Sigma)=\bbone\{\Pi=\Sigma\}
=\sum_{\Theta\in \Part(S)}\mu(\Pi:\Theta)\;\bbone\{\Theta\leq \Sigma\}
\end{equation}
holding for all $\Pi,\Sigma\in \Part(k)$. 
For $\CC$-valued random variables
$X_1,\dots,X_k$ with absolute moments of all orders it is well-known that 
\begin{equation}\label{equation:JointCumulantDef}
\kappa(X_1,\dots,X_k)=\sum_{\Pi\in \Part(k)}
\mu(\Pi:\one_k) \prod_{A\in \Pi}\Ebold\prod_{i\in A}X_i.
\end{equation}
This should be compared with formula \eqref{equation:JointCumulantDefPhys} above. 
\begin{Corollary}[The connected BKAR formula] \label{Corollary:BKAR}
In the setup of Theorem \ref{Theorem:BKAR} we have
\begin{equation}\label{equation:ConnectedBKAR}
\sum_{\Pi\in \Part(k)}\mu\left(\Pi:\one_k\right)f([\Pi])=
\sum_{\begin{subarray}{c}
\textup{\mbox{\scriptsize trees $\Tfrak$}}\\
\textup{\mbox{\scriptsize spanning $\langle k\rangle$}}
\end{subarray}}
\Ebold (\dd_\Tfrak f)(\weight_\Tfrak).
\end{equation}
\end{Corollary}
\proof Theorem \ref{Theorem:BKAR} has the trivial generalization
$$
f\left([\Pi]\right)
=\sum_{\begin{subarray}{c}
\mbox{\scriptsize forests $\Ffrak=(V,E)$ in $\langle k\rangle$}\\
\mbox{\scriptsize s.t. every edge $\{i,j\}\in E$ is}\\
\mbox{\scriptsize contained in some block of $\Pi$}
\end{subarray}}\Ebold(\dd_\Ffrak f)(\weight_\Ffrak)\;\;\;\mbox{for $\Pi\in \Part(k)$},
$$
whence the result  by M\"{o}bius inversion. \qed

\begin{Remark}
In the setup of Corollary \ref{Corollary:BKAR} it can be shown that
\begin{eqnarray*}
&&\sum_{\begin{subarray}{c}
\textup{\mbox{\scriptsize trees $\Tfrak$}}\\
\textup{\mbox{\scriptsize spanning $\langle k\rangle$}}
\end{subarray}}
\Ebold (\dd_\Tfrak f)(\weight_\Tfrak)\\
&=&\sum_{\begin{subarray}{c}
\mbox{\scriptsize maximal chains}\\
\one_k=\Phi_0>\cdots>\Phi_{k-1}=\zero_k\\
\mbox{\scriptsize in $\Part(k)$}
\end{subarray}}\int\cdots\int_{0=t_0<t_1<\cdots<t_{k-1}<t_k=1}\\
&&\frac{\partial^{k-1}}{\partial t_1\cdots \partial t_{k-1}}
f\left([\Phi_{k-1}]+\sum_{i=1}^{k-1}t_i([\Phi_{i-1}]-[\Phi_i])\right)\,\dd t_1\cdots \dd t_{k-1}.
\end{eqnarray*}
The proof is a straightforward if tedious application of Propositions \ref{Proposition:Obvious} and \ref{Proposition:ObviousConverse}, along with Remark \ref{Remark:CanBeShown}.
We note further that it is not too hard to verify \eqref{equation:ConnectedBKAR} directly by analyzing the Fundamental-Theorem-of-Calculus-induced relations standing between 
the integrals of the form appearing on the right side above and more general integrals of the same form attached to nonmaximal chains in $\Part(k)$.
\end{Remark}

\begin{Remark}
By a straightforward ``cutoff'' argument,
one can allow in equation \eqref{equation:ConnectedBKAR} any smooth  function $f$ defined on
an open set containing $\Sym_k^+$. We will use this freedom when proving Theorem \ref{Theorem:MalliavinApp}.
\end{Remark}
\subsection{Background on BKAR}  We provide pointers to the mathematical physics literature.
\subsubsection{Sample applications in physics}
While quite elementary, Theorem \ref{Theorem:BKAR} can be a very powerful tool when aptly deployed.
It already has a wide range of applications. Here is a non-exhaustive list of such:
\begin{itemize}
\item Explicit convergent expansions for the ground state energy of the spin-Boson model~\cite{AbdesselamSB}
as well as Nelson and polaron type models~\cite{AbdesselamH}.
\item Uniform $L^1$-type estimates for cumulants (truncated correlation functions)
in continuous spin models~\cite{AbdesselamPS}.
\item A close analog of the Ercolani-McLaughlin Theorem in the case of quartic interactions~\cite{Rivasseau}.
\item The construction of the L\'evy area for fractional Brownian motion with low Hurst exponent~\cite{MagnenU}.
\end{itemize}
In almost every application, the only difficulty is to ``see'' the quantity under study as
$f\left([\one_k]\right)$ for a suitable function $f$. 
\subsubsection{Cluster expansions} 
The typical application of Theorem \ref{Theorem:BKAR} is to provide an ``arboreal'' expression for 
the joint cumulant of
random variables $X_1,\ldots,X_k$ with a joint probability distribution (say) of the form $\exp(-\sum_{\alpha=1}^{k} V(x_\alpha)) \dd \mu_C(x)$
where $V$ is some function called a {\em potential} and $\dd \mu_C(x)$ is a Gaussian measure with covariance matrix $C$. Coefficients $s(i,j)$ are then introduced
as multipliers of the matrix elements $C(i,j)$
and can be viewed as the ``coupling'' between the random variables $X_i$ and $X_j$. The game is then to interpolate
between the fully coupled situation where $s(i,j)\equiv 1$
and the fully decoupled situation
where $s(i,j)\equiv \delta_{ij}$, i.e., the situation where the random variables become independent. 
The resulting decoupling expansion
generates partially decoupled intermediate terms corresponding to partitions $\Pi$ of $\langle k\rangle$ where random variables belonging
to different blocks of $\Pi$ become independent.  
This procedure is called a {\em cluster expansion} in the constructive quantum field theory
literature. The first such expansion appeared in~\cite{GlimmJS}. A simpler expansion was later introduced in~\cite{BrydgesF}.
Theorem \ref{Theorem:BKAR} belongs to a third generation of yet simpler cluster expansion formulas. 
\subsubsection{Generalizations of BKAR}
The utility of the representation of cumulants provided by the BKAR formula can be assessed by the quality of the bounds for such cumulants
that one can deduce from it. An important feature is that such expansions typically are in terms of trees instead of the much more numerous
connected graphs with arbitrarily many circuits. Indeed, BKAR produces a sum over minimally connected graphs. Actually Theorem \ref{Theorem:BKAR} is only  one of the simplest examples of a much more general hierarchy of combinatorial identities
which can generate hypergraphs instead of graphs and/or graphs which are (say)
minimally $p$-edge-connected. (The latter corresponds in quantum field theory jargon to  \linebreak $(p-1)$-particle irreducibility.) This very general framework was developed in~\cite{AbdesselamR2}.

\subsection{Proof of Theorem \ref{Theorem:MalliavinApp}}\label{subsection:ProofOfMalliavinApp}
We begin with a couple of lemmas, both of which are more or less standard.
\subsubsection{The functional $\widetilde{\HHH}$}
Given a tame function $f:\Mat_{k\times n}\rightarrow \CC$,
we define a function  $\widetilde{\HHH}f:\Sym^+_k\rightarrow \CC$ by the formula
$\widetilde{\HHH}f(Q)=\Ebold f(\sqrt{Q}Z)$
where as in the statement of Theorem \ref{Theorem:MalliavinApp}
the matrix $Z\in \Mat_{k\times n}$ is random with i.i.d.\ standard Gaussian entries.
Dominated convergence and continuity of the matrix square-root function implies that $\widetilde{\HHH}f$ is continuous.

\begin{Lemma}[The independent copies trick]\label{Lemma:IndependentCopies}
In the setup for Theorem \ref{Theorem:MalliavinApp} we have
$\kappa(f_1(\zeta),\dots,f_k(\zeta))=
\sum_{\Phi\in \Part(k)}\mu(\Phi:\one_k)(\widetilde{\HHH} f)([\Phi])$
where $f=f_1\otimes \cdots \otimes f_k$.
\end{Lemma}
\proof Let $\{\{\zeta_j^A\}_{j=1}^n\}_{\emptyset \neq A\subset \langle k\rangle}$ be an i.i.d.\ family of standard
Gaussian random variables. For $\Phi\in \Part(k)$ and $i\in \langle k\rangle$, let
$\Phi(i)\in \Phi$ be the block to which $i$ belongs. For $\Phi\in \Part(k)$ let $Z^\Phi\in \Mat_{k\times n}$
be the random matrix with entries $Z^\Phi(i,j)=\zeta^{\Phi(i)}_j$.
For each $\Phi\in \Part(k)$ the law of $[\Phi]^{1/2}Z$  is the same as that of $Z^\Phi$.
A straightforward calculation then shows that
$\Ebold(f_1\otimes \cdots \otimes f_k)([\Phi]^{1/2}Z)=
\prod_{A\in \Phi}\Ebold \prod_{i\in A} f_i(\zeta)$,
whence the result by \eqref{equation:JointCumulantDef} and the definition of $\widetilde{\HHH}$.
\qed

\subsubsection{The functional $\HHH$}
Let $\HHH f=\widetilde{\HHH}f\vert_{\Pos_k}$.
Dominated convergence shows that $\HHH f:\Pos_k\rightarrow\CC$ 
is once continuously differentiable, and the latter observation can be considerably
improved as follows.

\begin{Lemma}[The heat equation]\label{Lemma:HeatEquation}
For tame functions $f:\Mat_{k\times n}\rightarrow\CC$ and $i,i'\in \langle k\rangle$
we have $\sum_{j=1}^n\HHH  D_{ij}D_{i'j} f=
\nabla_{e_{ii'}+e_{i'i}}\HHH f$.
Consequently $\HHH f$ is infinitely differentiable and partial derivatives of $\HHH f$ of all orders extend continuously from $\Pos_k$ to $\Sym_k^+$.
\end{Lemma}
\proof This is proved by writing
$$\HHH f(Q)=\int_{\Mat_{k\times n}}f(X)
\frac{\exp(-\frac{1}{2}\trace \,X^\T Q^{-1}X)}{(2\pi)^{nk/2}(\det Q)^{n/2}}\prod_{i=1}^k\prod_{j=1}^n dX(i,j)$$
and taking advantage of the well-known system of partial differential equations
$$\left(
\underbrace{\nabla_{e_{ii'}+e_{i'i}}}_{\mbox{\scriptsize acts on entries of $Q$}}-\underbrace{\sum_{j=1}^nD_{ij}D_{i'j}}_{\mbox{\scriptsize acts on entries of $X$}}\right)
\frac{\exp(-\frac{1}{2}\trace\, X^\T Q^{-1}X)}{(2\pi)^{nk/2}(\det Q)^{n/2}}\equiv 0.$$
 We can safely omit further details.
 \qed
\subsubsection{Concluding calculation}
In the setup for Theorem \ref{Theorem:MalliavinApp},
where we write \linebreak $f=f_1\otimes \cdots \otimes f_k$ and let $\sum_\Tfrak$ denote
summation over all trees $\Tfrak$ spanning $\langle k\rangle$, we have
\begin{eqnarray*}
&&\sum_\Tfrak\Ebold\, (L_\Tfrak f)(\sqrt{\weight_\Tfrak}\,Z)\,=\,\sum_\Tfrak\Ebold(\HHH L_\Tfrak f)(\weight_\Tfrak)\,
=\,\sum_\Tfrak\Ebold (\dd_\Tfrak\HHH f)(\weight_T)\\
&=&\lim_{\epsilon\downarrow 0}\sum_\Tfrak\Ebold (\dd_\Tfrak\HHH f)(\weight_T+\epsilon \Ibold_k)\,=\,\lim_{\epsilon\downarrow 0}\sum_{\Phi\in \Part(k)}\mu(\Phi:\one_k)(\HHH f)([\Phi]+\epsilon\Ibold_k)\\
&=&\sum_{\Phi\in \Part(k)}\mu(\Phi:\one_k)(\widetilde{\HHH}f)([\Phi])\,=\,\kappa(f_1(\zeta),\dots,f_k(\zeta)).
\end{eqnarray*}
The first step is justified by Fubini's theorem, the second by Lemma \ref{Lemma:HeatEquation},
the third by Lemma \ref{Lemma:HeatEquation} and dominated convergence,
the fourth by connected BKAR (Corollary \ref{Corollary:BKAR}), the fifth by continuity of $\widetilde{\HHH}f$,
and the last by Lemma \ref{Lemma:IndependentCopies}.
\qed

\section{Proof of  the main result}
\label{section:OperatorTheoreticCharacterization}

\subsection{The splicing identity}\label{subsection:SplicingIdentity}
We execute the derivative calculation needed to apply Theorem \ref{Theorem:MalliavinApp}
to the analysis of the expression under the limit sign on the  right side of \eqref{equation:MatricialFormOftHooft}. 
\subsubsection{Setup}
We fix $\theta\in \Perm(n)$, $\gamma\in \Color(n)$
and an onto $\theta$-invariant map $\nu:\langle n\rangle\rightarrow\langle k\rangle$
as in the setup for Theorem \ref{Theorem:BigKahuna}. We put $k=c(\theta)$.
In addition we  fix a tree $\Tfrak$ spanning $\langle k\rangle$ and a positive integer $N$. 
As in the statement of Theorem~\ref{Theorem:tHooft} we fix a decomposition
$$\theta=(i_{1,1}\cdots i_{1,n_1})\cdots (i_{k,1}\cdots i_{k,n_k})\;\;\;(n_1+\cdots+n_k=n)$$
into cycles. After relabeling if necessary, we may and we do assume that 
$\nu(i_{\alpha,\beta})=\alpha$
for $(\alpha,\beta)$ such that $1\leq \alpha \leq k$ and $1\leq \beta\leq n_\alpha$.

\subsubsection{Commutative variables and a differential operator}
Let
$$\{\{\{x_{ij}(\alpha,\beta)\}_{\alpha,\beta=1}^N\}_{i=1}^k\}_{j=1}^n$$
be a family of independent commutative algebraic variables.
Consider also the differential operator
$$
L_\Tfrak^{(N)}=\prod_{\begin{subarray}{c}
\textup{\mbox{\scriptsize edges $\{i,i'\}$}}\\
\textup{\mbox{\scriptsize of $\Tfrak$}}
\end{subarray}}\;\;\sum_{j=1}^n\sum_{\alpha,\beta=1}^N\frac{\partial^2}{\partial x_{ij}(\alpha,\beta)
\partial x_{i'j}(\alpha,\beta)}
$$
acting on the commutative polynomial algebra generated by the variables $x_{ij}(\alpha,\beta)$
over the complex numbers.
\subsubsection{Matrices with commutative polynomial entries}
Put
$$\hat{x}_{ij}(\alpha,\beta)=
\left\{\begin{array}{rl}
x_{ij}(\alpha,\alpha)&\mbox{if $\alpha=\beta$,}\\
\frac{x_{ij}(\alpha,\beta)+\ii x_{ij}(\beta,\alpha)}{\sqrt{2}}&\mbox{if $\alpha<\beta$,}\\
\frac{x_{ij}(\beta,\alpha)-\ii x_{ij}(\alpha,\beta)}{\sqrt{2}}&\mbox{if $\alpha>\beta$}
\end{array}\right.
$$
and let 
$$X^{(N)}_{ij}=\left[\begin{array}{cccc}
\hat{x}_{ij}(1,1)&\dots&\hat{x}_{ij}(1,N)\\
\vdots&&\vdots\\
\hat{x}_{ij}(N,1)&\dots&\hat{x}_{ij}(N,N)
\end{array}\right]\;\;\mbox{for $i=1,\dots,k$ and $j=1,\dots,n$.}$$
We remark that if one evaluates the algebraic variables $x_{ij}(\alpha,\beta)$ at i.i.d.\ standard
Gaussian random variables,
then the matrices $X^{(N)}_{ij}$ evaluate to independent \linebreak $N$-by-$N$ GUE matrices. 
However, we do not yet carry out that evaluation. We are concerned for the moment only
with differentiation of polynomial expressions.

\begin{Proposition}\label{Proposition:Ghastly}
Notation and assumptions are as above. We have a relation
\begin{eqnarray}\label{equation:Ghastly}
&&L_\Tfrak^{(N)}\left(\trace (X^{(N)}_{1,\gamma(i_{1,1})}\cdots X^{(N)}_{1,\gamma(i_{1,n_1})})\cdots
\trace (X^{(N)}_{k,\gamma(i_{k,1})}\cdots X^{(N)}_{k,\gamma(i_{k,n_k})})\right)\\
\nonumber&=&\trace\,\left(
\Poly_{\theta,\gamma,\nu,\Tfrak}\bigg\vert_{\textup{\mbox{\scriptsize $\Xbold_i=X_{\nu(i),\gamma(i)}^{(N)}$
for $i\in \langle n\rangle$}}}\right).
\end{eqnarray}
\end{Proposition}
\noindent This proposition is the sole motivation for the 
definition of $\Poly_{\theta,\gamma,\nu,\Tfrak}$.
\proof We work in a setup similar to that of the proof of Proposition \ref{Proposition:SpliceCount}.
Let 
\begin{eqnarray*}
\vec{\Tfrak}&=&\{(i,i')\in \langle k\rangle^2\mid \{i,i'\}:\;\mbox{edge of $\Tfrak$}\},\\
\pi&=&((i,i')\mapsto i):\vec{\Tfrak}\rightarrow\langle k\rangle,\\
\rho&=&((i,i')\mapsto (i',i)):\vec{\Tfrak}\rightarrow\vec{\Tfrak}.
\end{eqnarray*}
After  opening the brackets in the definition  of $L_\Tfrak^{(N)}$
and noting that
$$\frac{\partial^2}{\partial x_{ij}(\alpha,\beta)\partial x_{i'j}(\alpha,\beta)}=
\frac{\partial^2}{\partial \hat{x}_{ij}(\alpha,\beta)\partial \hat{x}_{i'j}(\beta,\alpha)},$$
 we get the formula
$$L_\Tfrak^{(N)}=
\sum_{
\kbold:\vec{\Tfrak}\rightarrow\langle N\rangle}\;\;
\sum_{\begin{subarray}{c}
\jbold:\vec{\Tfrak}\rightarrow \langle n\rangle\\
\st\,\jbold\circ \rho=\jbold
\end{subarray}}\;\;
\prod_{\vec{e}\in \vec{\Tfrak}}\;
\frac{\partial}{\partial \hat{x}_{\pi(\vec{e}),\jbold(\vec{e})}(\kbold(\vec{e}),\kbold(\rho(\vec{e})))}.
$$
By similarly opening the brackets we obtain identities
\begin{eqnarray*}
&&\trace (X^{(N)}_{1,\gamma(i_{1,1})}\cdots X^{(N)}_{1,\gamma(i_{1,n_1})})\cdots
\trace (X^{(N)}_{k,\gamma(i_{k,1})}\cdots X^{(N)}_{k,\gamma(i_{k,n_k})})\\
&=&\sum_{\ibold:\langle n\rangle\rightarrow \langle N\rangle}\;\;\prod_{s\in \langle n\rangle}\;
\hat{x}_{\nu(s),\gamma(s)}(\ibold(s),\ibold(\theta(s))),\\\\\\
&&\trace\,\left(
\Poly_{\theta,\gamma,\nu,\Tfrak}\bigg\vert_{\textup{\mbox{\scriptsize $\Xbold_i=X_{\nu(i),\gamma(i)}^{(N)}$
for $i\in \langle n\rangle$}}}\right)\\
&=&\sum_{
\tau\in \Splice_\Tfrak(\nu,\gamma)}\;\;
\sum_{
\ibold:\langle n\rangle\rightarrow \langle N\rangle}\;\;\prod_{s\in \langle n\rangle}
\left\{\begin{array}{rl}
\hat{x}_{\nu(s),\gamma(s)}(\ibold(s),\ibold(\theta\tau(s)))&\mbox{if $\tau(s)=s$,}\\
\bbone\left\{\ibold(s)=\ibold(\theta\tau(s))\right\}&\mbox{if $\tau(s)\neq s$.}
\end{array}\right.
\end{eqnarray*}
We thus have
\begin{eqnarray*}
&&(\mbox{LHS of \eqref{equation:Ghastly}})\\
&=&\sum_{\ibold:\langle n\rangle\rightarrow \langle N\rangle}\sum_{
\kbold:\vec{\Tfrak}\rightarrow\langle N\rangle}\;\;
\sum_{\begin{subarray}{c}
\jbold:\vec{\Tfrak}\rightarrow \langle n\rangle\\
\st\,\jbold\circ \rho=\jbold
\end{subarray}}\;\;\sum_{\begin{subarray}{c}
\psi:\vec{\Tfrak}\rightarrow\langle n\rangle\\
\mbox{\scriptsize one-to-one}
\end{subarray}}\\
&&
\qquad \bbone\{\pi=\nu\circ \psi\}\bbone\{\jbold=\gamma\circ \psi\}
\bbone\{\kbold=\ibold\circ \psi\}\bbone\{\kbold\circ \rho=\ibold\circ \theta\circ\psi\}\\
& & \qquad
\prod_{s\in \langle n\rangle\setminus \psi(\vec{\Tfrak})}\hat{x}_{\nu(s),\gamma(s)}(\ibold(s),\ibold(\theta(s)))\\
&=&
\sum_{\begin{subarray}{c}
\psi:\vec{\Tfrak}\rightarrow\langle n\rangle\\
\mbox{\scriptsize one-to-one}\\
\st \;\nu\circ \psi=\pi\\
\mbox{\scriptsize and}\;\gamma\circ \psi\circ \rho=\gamma \circ \psi
\end{subarray}}\;\;\;
\sum_{\begin{subarray}{c}
\ibold:\langle n\rangle\rightarrow \langle N\rangle\\
\st\,\ibold\circ \psi\circ \rho=\ibold\circ \theta\circ \psi
\end{subarray}}\;\;\prod_{s\in \langle n\rangle\setminus \psi(\vec{\Tfrak})}\hat{x}_{\nu(s),\gamma(s)}(\ibold(s),\ibold(\theta(s)))\\
&=&\sum_{
\tau\in \Splice_\Tfrak(\nu,\gamma)}\;\;
\sum_{\begin{subarray}{c}
\ibold:\langle n\rangle\rightarrow \langle N\rangle\\
\st\,\ibold(\theta\tau(s))=\ibold(s)\\
\textup{\mbox{\scriptsize for $s\in \langle n\rangle$ s.t. $\tau(s)\neq s$}}
\end{subarray}}\;\;\prod_{\begin{subarray}{c}
s\in \langle n\rangle\\
\mbox{\scriptsize s.t.}\;\tau(s)=s
\end{subarray}}
\hat{x}_{\nu(s),\gamma(s)}(\ibold(s),\ibold(\theta(s)))\\
&=&(\mbox{RHS of \eqref{equation:Ghastly}}),
\end{eqnarray*}
which finishes the proof. \qed

\subsection{Concluding arguments}

\subsubsection{Rewrite of the right side of \eqref{equation:MatricialFormOftHooft}}
Let $\Xi^{(N)}_{ij}$ be the $N$-by-$N$ GUE matrix to which the matrix $X^{(N)}_{ij}$
evaluates when all the variables $x_{ij}(\alpha,\beta)$ are evaluated at i.i.d. standard normal random variables.
Then by Theorem \ref{Theorem:MalliavinApp} combined with Proposition \ref{Proposition:Ghastly},
with $\sum_\Tfrak$ denoting summation over trees $\Tfrak$ spanning $\langle k\rangle$
and picking the random matrices $\weight_\Tfrak$
independent of all the GUE matrices $\Xi_{ij}^{(N)}$, 
we have the formula
\begin{eqnarray}\label{equation:ExactAndScary}
&&\left(\begin{array}{l}
\mbox{expression under the limit on the right}\\
\mbox{side of \eqref{equation:MatricialFormOftHooft}
multiplied by $N^{2+n/2-k}$}\end{array}\right)\\
\nonumber&=&
\sum_\Tfrak \Ebold \,\trace\left(
\left(\Poly_{\theta,\gamma,\nu,\Tfrak}\left|_{
\begin{subarray}{l}
\Xbold_i=\Zbold_\Tfrak(\nu(i),\gamma(i))\\
\mbox{\scriptsize for $i\in \langle n\rangle$}
\end{subarray}}\right)\right.\left|_{\begin{subarray}{l}\Zbold(i,j)=\Xi_{ij}^{(N)}\\
\mbox{\scriptsize for $(i,j)\in \langle k\rangle\!\times\! \langle n\rangle$}
\end{subarray}}\right)\right.\\
\nonumber&=&\Ebold\, \trace\left( \Arb_{\theta,\gamma,\nu}\left|_{\begin{subarray}{l}\Zbold(i,j)=\Xi_{ij}^{(N)}\\
\mbox{\scriptsize for $(i,j)\in \langle k\rangle\!\times\! \langle n\rangle$}
\end{subarray}}\right).\right.
\end{eqnarray}

\subsubsection{Passage to the limit}
Divide through \eqref{equation:ExactAndScary} by $N^{2+n/2-k}$.
Consider the limit as $N\rightarrow\infty$ on extreme left and extreme right.
On the one hand, by Theorem \ref{Theorem:tHooft} on the left, we recover the quantity $|\Map_0(\theta,\gamma)|$.
On the other hand, by Voiculescu's result \cite[Thm. 2.2]{Voiculescu} on the right,
we recover the right side of \eqref{equation:BigKahuna}.
The proof of Theorem \ref{Theorem:BigKahuna}  is now complete.
  \qed

\begin{Remark}
The preceding argument provides a proof of the existence of the limit
on the right side of \eqref{equation:MatricialFormOftHooft} independent of the usual diagram-based 
approach.
\end{Remark}
\begin{Remark}\label{Remark:Susceptible}
Since relation \eqref{equation:ExactAndScary} already holds exactly for each fixed $N$,
the methods introduced in this paper seem to be susceptible to higher genus generalization.
\end{Remark}


\end{document}